\documentclass[12pt,a4paper]{amsart}
\usepackage{xypic,amssymb,amsthm,a4wide}
\def\F{{\mathcal F}}\def\O{{\mathcal O}}\def\Pc{{\mathcal 
P}}\def\M{{\mathcal M}}
\def\Qc{{\mathcal Q}}
\def\A{{\mathcal A}}\def\E{{\mathcal E}}\def\G{{\mathcal G}}
\def\Nc{{\mathcal N}}\def\L{{\mathcal L}}\def\K{{\mathcal K}}

\def\bS{{\mathbf S}}
\newtheorem{theorem}{Theorem}[section]
\newtheorem{corollary}[theorem]{Corollary}
\newtheorem{lemma}[theorem]{Lemma}
\newtheorem{proposition}[theorem]{Proposition}
\newtheorem{definition}[theorem]{Definition}
\newtheorem{rem}[theorem]{Remark}
\newenvironment{remark}{\begin{rem}\rm}
       {\hfill$\vartriangle$\end{rem}}
\begin{document}
\thispagestyle{empty}
\addvspace{20pt}
\begin{center}{\bfseries\large RELATIVELY STABLE BUNDLES OVER \\[10pt] ELLIPTIC FIBRATIONS}
\\[30pt]  {\sc Claudio Bartocci\thinspace \ddag \\[3pt]  Ugo Bruzzo\thinspace \S$\star$ \\[3pt] Daniel
Hern\'andez Ruip\'erez , Jos\'e M.~Mu\~noz Porras\thinspace \P\\[10pt]}
 \ddag\thinspace Dipartimento di Matematica, Universit\`a di
Ge\-no\-va, \\ Via Dodecaneso 35, 16146 Ge\-no\-va, Italy \\[3pt]
\S\thinspace Scuola Internazionale Superiore di Studi Avanzati, \\  Via
Beirut 2-4, 34014 Trieste, Italy \\[3pt] $\star$\thinspace Istituto Nazionale di Fisica Nucleare, \\
Sezione di Trieste\\[3pt]
 \P\thinspace Departamento de Matem\'aticas,
Universidad de Salamanca, \\ Plaza de la Merced 1-4, 37008 Salamanca, Spain
\\[10pt]  {\tt bartocci@dima.unige.it, bruzzo@sissa.it,\\ ruiperez@usal.es, jmp@usal.es}
\end{center}
\vfill\begin{center}
20 February 2000, Revised 25 August 2001 and 13 September 2004
\end{center}
\vfill
\begin{quote} {\sc Abstract.}
We consider a relative Fourier-Mukai transform defined on  elliptic fibrations
over an arbitrary base scheme. This is used to construct relative 
Atiyah sheaves and generalize Atiyah's and Tu's results about semistable
sheaves over elliptic curves to the case of elliptic fibrations. Moreover we
show that this transform preserves relative (semi)stability of sheaves of
positive relative degree.
\end{quote}
\vfill{\small
This research was partially supported by the Spanish DGES through the research
project PB96-1305, by ``Junta de Castilla y Le\'on''  through the 
research project
SA27/98, by GNSAGA (Italian National Research Council) and by an 
Italian-Spanish
cooperation project. The authors are  members of VBAC (Vector Bundles 
on Algebraic
Curves), which is partially supported by {\sc Eager} (EC FP5 Contract no.
HPRN-CT-2000-00099) and  {\sc Edge} (EC FP5 Contract no. HPRN-CT-2000-00101).}
\newpage
\section{Introduction} In \cite{BBHM} we defined a relative Fourier-Mukai
transform on elliptic K3 surfaces. A relative Fourier-Mukai transform 
for abelian
schemes had been considered by Mukai \cite{Muk2} generalizing  his previous
definition holding for abelian varieties \cite{Muk1}. The case of (absolute)
transforms for K3 surfaces can be found in
\cite{BBH} and \cite{Muk3}.

In \cite{HM} a relative Fourier-Mukai transform was introduced on elliptic
fibrations
$p\colon X\to B$ over an arbitrary normal base scheme $B$ (see also 
\cite{Br,BM}
for the case of elliptic surfaces and that of elliptic fibrations over smooth
bases). In this paper we prove that this transform preserves the relative
stability of the sheaves on
$X$ of positive relative degree. This, together with the results in \cite{HM},
generalizes to elliptic fibrations the classical theory of moduli 
spaces of vector
bundles on elliptic curves \cite{A,Tu}.

As it is well known, every indecomposable vector bundle over a smooth 
elliptic curve $C$
is semistable,  and stability holds if and only if   rank and degree are
relatively prime. This result is proved in
\cite{Tu} by showing that the Atiyah bundle $\mathcal A_n$ on $C$ is 
semistable.
The Atiyah bundle
$\mathcal A_n$ is the rank $n$ indecomposable vector bundle defined 
inductively by
the exact sequence $$0
\to \mathcal O_C \to \mathcal A_n \to \mathcal A_{n-1} \to 0$$ where 
$\mathcal A_0
=\mathcal O_C$. These constructions carry over to the case of elliptic curves
with a singular point. As a direct consequence of this characterization,  the 
moduli space of semistable bundles of rank $n$ and degree $d$ is naturally
isomorphic to the symmetric product $S^h C$, where $h$ is the g.c.d. of $n$ and
$d$.

It turns out that these results can be extended to the case of 
elliptic fibrations
as it is shown in
\cite{HM} and in the present paper (see also \cite{Br,BM}). In our construction
the relative Fourier-Mukai transform plays a central role. In particular, this
transform allows us to give a natural definition of the relative rank 
$i$ Atiyah
sheaf $\mathcal A_i$    (Prop.
\ref{atiyah} and Def. \ref{relatiyah}) and to obtain explicitly  the 
isomorphism
between the compactified relative Jacobian $\bar J_n$ of (flat families of)
torsion-free, rank $1$ and degree $n$ sheaves and the moduli space $\mathcal M
(n,-1)$ of rank $n$ relatively stable sheaves of degree $-1$ (Cor. 
\ref{cor3.5}).

Moreover, the results in this paper provide a kind of duality between the
following functors.

(1) The first functor will be denoted by ${\mathbf F}^i_{X}(n,d)$; it 
associates
to every scheme $S$ over $B$ the scheme paramet\-rizing exact 
sequences of sheaves
on $X_S$
\begin{equation} 0\to\L\otimes\A_i \to\E_n \to\E_{n-i}\to 0 \label{flag}
\end{equation} where $\L$ is torsion-free, of rank 1 and relative 
degree zero and
$\E_n$,
$\E_{n-i}$ are relatively stable sheaves of relative degree $d\ge 0$ 
and rank $n$
and $n-i$, respectively.

(2) The second functor, denoted by ${\mathbf H}^i_{X}(r,d)$, 
associates to $S$ the
scheme param\-etrizing the exact sequences of sheaves on $X_S$ \begin{equation}
0\to\E_r\to\E'_r\to\K(x)\to 0\,,
\label{higher}
\end{equation} where $\E_r$, $\E'_r$ are relatively stable sheaves of 
rank $r$ and
$\K(x)$ has length $i$ on every fibre and is concentrated on the 
image $x(S)$ of a
section $x\colon S\hookrightarrow X_S$ of $p_S$; $d$ denotes the 
relative degree
of $\E_r$.

In Theorem \ref{main} we prove that the Fourier-Mukai transform establishes an
isomorphism between the functors ${\mathbf F}^i_{X}(n,d)$ and ${\mathbf
H}^i_{{\widehat X}}(d,-n)$, where ${\widehat X}$ is the compactified relative
Jacobian of $X$.

\medskip
\textbf{Acknowledgement.}
We thank Igor Burban and Bernd Kreussler for pointing out that in a previous version of this paper Proposition \ref{atiyah} was incorrectly stated.

\section{Preliminaries\label{prelim}} Let $p\colon X\to B$ be an elliptic
fibration, where $B$ is a normal algebraic scheme over an algebraically closed
field, and the fibres of the projection $p$ are geometrically 
integral Gorenstein
curves of arithmetic genus 1. We also assume that $p$ has a 
section $e$ taking
values in the smooth locus of $p$. We write
$H=e(B)$; this is a relative polarization. We shall also denote by 
$X_t$ the fibre
of $p$ over $t\in B$.

The sheaf $\omega=R^1p_\ast\O_X$ is a line bundle on $B$. By relative duality,
$\omega\simeq (p_\ast \omega_{X/B})^\ast$, where $\omega_{X/B}$ is the relative
dualizing sheaf, and also
$\omega_{X/B}\simeq p^\ast\omega^{-1}$.

We denote by $\widehat X$ the compactified relative Jacobian of $X$. 
This is the
scheme representing the functor which to any scheme morphism $f\colon S\to B$
associates the space of equivalence classes of $S$-flat sheaves on $p_S\colon
X\times_B S\to S$, whose restrictions to the fibres of
$f$ are torsion-free, of rank one and degree zero; two such sheaves 
$\F$, $\F'$ are
considered to be equivalent if $\F'\simeq\F\otimes p_S^\ast\Nc$ for a 
line bundle
$\Nc$ on
$S$ (cf. \cite{AK}).

There is a natural isomorphism of $B$-schemes
\begin{eqnarray*} \varpi\colon\ X & \to & {\widehat X} \\ x & \mapsto 
& {\mathfrak
m}_x^\ast\otimes\O_{X_t}(-e(t)) \end{eqnarray*} where ${\mathfrak m}_x$ is the
ideal sheaf of the point $x$ in $X_t$. The section $e$ induces then a section
$\hat e$ of $\hat p\colon\widehat X\to B$ whose image is a relative 
polarization
$\Theta$.

We can normalize the Poincar\'e sheaf $\Pc$ on $X\times_B{\widehat 
X}$ so that $$
\Pc_{\vert{H\times_B {\widehat X}}}\simeq \O_{{\widehat X}}\,. $$ One has
$$
\Pc\simeq\mathcal J\otimes \pi^\ast\O_X(H)\otimes\hat\pi^\ast\O_{{\widehat
X}}(\Theta)\otimes q^\ast\omega^{-1}\,,
$$ where $\mathcal J$ is the ideal of the graph $\gamma\colon X\hookrightarrow
X\times_B {\widehat  X}$ of
$\varpi\colon X\to{\widehat  X}$ and $q=p\circ\pi=\hat 
p\circ\hat\pi$. Then  $\Pc$
is flat over both $X$ an $\widehat X$. Moreover, one has:

\begin{lemma} The dual $\Pc^\ast$ of the Poincar\'e sheaf is flat over both
factors and
$\Pc$ is reflexive, that is, $\Pc\simeq\Pc^{\ast\ast}$. \end{lemma}
\begin{proof} One has ${\mathcal
E}xt\,^1_{\O_{X_t}}(\Pc_\xi,\O_{X_t})=0$  for every point $\xi\in \widehat X$
(where $t=\tilde p(\xi)$), as it was proved in Lemma 0.4 of \cite{FMW}.
Now, by Theorem 1.10 of \cite{AK}, ${\mathcal
E}xt\,^1_{\O_{X\times_B  {\widehat
X}}}(\Pc,\O_{X\times_B {\widehat  X}})=0$, and the base change property for the
local Ext's (\cite{AK}, Theorem 1.9) implies that
$(\Pc^\ast)_\xi\simeq (\Pc_\xi)^\ast$ for every point $\xi\in{\widehat  X}$. By
applying once more Theorem 1.9 of \cite{AK}, we have that $\Pc^\ast$ is flat
over ${\widehat  X}$.  After interchanging the roles of $X$ and ${\widehat  X}$,
we obtain that $\Pc^\ast$ is flat over
$X$ as well. Then, $\Pc^\ast$ defines a morphism $\iota\colon{\widehat  X}\to
{\widehat  X}$ that maps any rank-one torsion-free  zero-degree coherent sheaf
$\F$ on a fibre
$X_t$ to its dual $\F^\ast$. Thus,
$(1\times\iota)^\ast\Pc\simeq\Pc^\ast\otimes\hat\pi^\ast\Nc$ for some 
line bundle
$\Nc$ on
${\widehat  X}$, which turns out to be trivial. Then \begin{equation}
(1\times\iota)^\ast\Pc\simeq\Pc^\ast\,.
\label{e:inv}
\end{equation} The morphism
$\iota\circ\iota\colon{\widehat  X}\to{\widehat  X}$ is the identity on the
relative Jacobian
$J(X/B)\subset{\widehat  X}$; by separateness 
$\iota\circ\iota=\mbox{Id}$, and the
isomorphism (\ref{e:inv}) implies
$\Pc\simeq\Pc^{\ast\ast}\,.$
\end{proof}

As a consequence, every coherent sheaf $\F$ on $X\times_B  S$ flat 
over $S$ whose
restrictions to the fibres of
$X\times_B  S \to S$ are torsion-free and of rank one and degree zero 
is reflexive,
$\F\simeq\F^{\ast\ast}$.

For any morphism $S\to B$ of algebraic varieties we define the (relative)
Fourier-Mukai transform as the functor $\bS_S\colon D^-(X_S)\to 
D^-(\widehat X_S)$
(where $D^-$ denotes the derived category of coherent sheaves consisting of
complexes bounded from above, and a subscript $S$ denotes base change 
to $S$, so
that
$X_S=X\times_B S$ etc.) given by
$$\bS_S(F)={\mathbf R}\hat\pi_{S\ast}(\pi_S^\ast F \ \stackrel{L}{\otimes}\
\Pc_S)\,,$$
$\pi\colon X\times_B\widehat X\to X$ and $\hat \pi\colon X\times_B\widehat X\to
\widehat X$ being the natural projections. We can also define the Fourier-Mukai
transforms of a single sheaf $\mathcal F$ as
$$
\bS_S^i(\mathcal F)=R^i\hat\pi_{S\ast}(\pi_S^\ast \mathcal F \otimes 
\Pc_S) $$ One
then says that a sheaf
$\mathcal F$ on
$X_S$ is WIT$_i$ if
$$\bS_S^j(\mathcal F)=0\quad\mbox{for}\ j\ne i\,.$$ In this case we write
$\widehat \F =
\bS_S^i(\mathcal F)$.

Several results that were proved in \cite{BBHM} in the case of elliptic K3
surfaces also hold under the conditions of the present paper. For the reader's
convenience we recall here some of them; further details can be found in
\cite{BBHM}.

\begin{proposition} Let $\F$ be a sheaf on $X_S$, flat over $S$. The 
Fourier-Mukai
transforms
$\bS_S^i (\F)$ are $S$-flat as well. Moreover, for every morphism 
$g\colon T\to S$
one has
$g_{{\widehat X}}^\ast\bS_S^1(\F)\simeq \bS_T^1(g_X^\ast\F)$, where
$g_X\colon X_T\to X_S$, $g_{{\widehat X}}\colon{\widehat X}_T\to{\widehat X}_S$
are the  morphisms induced by $g$.
\label{basechange}
\qed\end{proposition}

Let $f\colon S\to B$ be a morphism and $\L$ a coherent sheaf on 
$X\times_B S$ flat
over
$S$ whose restrictions to the fibres of $p_S$ are torsion-free and 
have rank one
and degree zero. Let
$\phi\colon S\to{\widehat X}$ be the morphism determined by the universal
property, so that
$$ (1\times\phi)^\ast\Pc\simeq\L\otimes p_S^\ast\Nc$$ for a line bundle
$\Nc$ on $S$. Let
$\Gamma\colon S\hookrightarrow {\widehat X}_S$ be the graph of the morphism
$\iota\circ\phi\colon S\to{\widehat X}$, where $\iota$ is the natural 
involution of
${\widehat X}$ induced by the operation of taking the dual.

Lemma 2.11 of \cite{BBHM} takes now the form: \begin{lemma} In the 
above situation
$\bS^0_S(\L)=0$ and
$\bS_S^1(\L)\otimes\hat p_S^\ast\Nc\simeq \Gamma_{\ast}(f^\ast\omega)$. In
particular,
\begin{list}{}{\itemsep=2pt}
\item[(1)] $\bS_{{\widehat X}}^0(\Pc)=0$ and $\bS_{{\widehat X}}^1(\Pc)\simeq
\zeta_\ast\hat p^\ast\omega$, where $\zeta\colon {\widehat X}\hookrightarrow
{\widehat X}\times_B{\widehat X}$ is the graph of the morphism 
$\iota$. \item[(2)]
$\bS_{{\widehat X}}^0(\Pc^\ast)=0$ and
$\bS_{{\widehat X}}^1(\Pc^\ast)\simeq \delta_\ast\hat p^\ast\omega$, where
$\delta\colon{\widehat X}\hookrightarrow{\widehat X}\times_B{\widehat 
X}$ is the
diagonal immersion.
\end{list}
\label{SSNew}
\qed\end{lemma} We also recover Corollary 2.12 of \cite{BBHM}: 
\begin{corollary}
Let $\L$ be a rank-one,  zero-degree, torsion-free coherent sheaf on a fibre
$X_t$. Then
$$
\bS^0_t(\L)=0\,,\qquad \bS^1_t(\L)=\kappa([\L^\ast])\,, $$ where $[\L^\ast]$ is
the point of
${\widehat X}_t$ defined by $\L^\ast$. \qed\label{cor1}\end{corollary}

By reversing the roles of $X$ and $\widehat X$, and using the sheaf
$\Qc=\Pc^\ast\otimes\pi^\ast p^\ast\omega^{-1}$ instead  of $\Pc$, we define a
functor\footnote{The present definition of the functor
$\hat\bS$ is slightly different than in
\cite{BBHM} due to the presence of the factor $p^\ast\omega^{-1}$.}
\begin{eqnarray*}
\widehat\bS_S\colon D^-(\widehat X_S)&\to& D^-(X_S) \\ G&\mapsto&
R\pi_{S\ast}(\hat\pi_S^\ast G\stackrel{L}{\otimes}\Qc_S) 
\end{eqnarray*} and the
corresponding functors $\widehat\bS_S^i$, $i=0,1$. Proceeding as in 
Theorem 3.2 of
\cite{BBHM},  and taking into account Lemma
\ref{SSNew}, we obtain an invertibility result: \begin{proposition}  For every
$G\in D^-({\widehat X}_S)$,
$F\in D^-(X_S)$ there are functorial isomorphisms
$$
\bS_S (\widehat{\bS}_S(G))\simeq G[-1]\,,\quad \widehat{\bS}_S(\bS_S (F))\simeq
F[-1]
$$ in the derived categories $D^-({\widehat X}_S)$ and $D^-(X_S)$, 
respectively.
\label{invert}
\qed\end{proposition}

\begin{proposition} Let $\F$ be a rank $n$ coherent sheaf on $X_S$, 
flat over $S$,
of relative degree $d$.  The relative Chern character of the Fourier-Mukai
transform
$\bS_S(\F)$ is $(d,-n)$. \label{chern}
\end{proposition}
\begin{proof} To compute the relative invariants of $\bS_S(\F)$ we 
take for $S$ a
point
$t\in B$.  The morphism
$X_t\times{\widehat X}_t\to {\widehat X}_t$ is a local complete intersection
morphism in the sense of \cite{Fu}, so that one can  apply the Riemann-Roch
theorem for singular varieties applies (Corollary 18.3.1 of
\cite{Fu}) . \end{proof} We shall denote by $d(\F)$ the relative 
degree of a sheaf
$\F$ and by $\mu(\F)=d(\F)/\mbox{rk}\,(\F)$ its  relative slope.

\begin{corollary} If $\F$ is WIT$_i$ and $d\ne 0$, then
$\mu(\widehat\F)=-1/\mu(\F)$.\label{mu}
\end{corollary}
\begin{corollary}
\begin{enumerate}
\item[(1)] If $\F$ is WIT$_0$, then $d(\F)\ge 0$, and
$d(\F)=0$ if and only if $\F=0$.
\item[(2)] If $\F$ is WIT$_1$, then $d(\F)\le 0$. \end{enumerate}
\label{wit}
\end{corollary}

Proceeding as in (\cite{A}, Theorem 5,  or
\cite{Tu}) (cf.~also   \cite{FMW}, Lemma 3.1), one proves the following
Proposition, which  allows one to define the rank
$n$ Atiyah  sheaf on a fibre
$X_t$.
\begin{proposition} For every $n\ge0$, there exists a unique
torsion-free indecomposable semistable rank
$n$ and degree 0 sheaf $\A_n^{(t)}$ on $X_t$ which satisfies
$H^0(X_t,\A_n^{(t)})\ne 0$. Moreover  there is an exact sequence
$$ 0 \to \O_{X_t} \to \A_n^{(t)} \to \A_{n-1}^{(t)} \to 0\,.
$$
\label{atiyah}
\qed\end{proposition}
\begin{corollary} The Atiyah sheaf
$\A_n^{(t)}$ is WIT$_1$ and $\bS_t^1(\A_n^{(t)})=\O_{n\hat e(t)}$.
\end{corollary}
\begin{proof} One applies Corollary \ref{cor1} to the exact sequence 
in Proposition
\ref{atiyah} and  uses induction. \end{proof} Then, by the invertibility of the
Fourier-Mukai transform, we have
$\A_n^{(t)}=\hat\bS_t^0(\O_{n\hat e(t)})$. This motivates the following:
\begin{definition} The relative rank $n$ Atiyah sheaf associated with the
fibration $p\colon X$ $\to B$ is the sheaf $$
\A_n=\hat\bS_B^0(\O_{n\Theta})\otimes p^\ast\omega\,. $$ \label{relatiyah}
\end{definition}
\noindent The factor $p^\ast\omega$ is introduced so that $\A_1\simeq\O_X$.
\section{Preservation of stability} In this section we prove that the
Fourier-Mukai transform preserves the (semi)sta\-bil\-ity of sheaves 
of positive
relative degree. The zero degree case is somehow more delicate and  has been
treated separately in \cite{HM}.

When considering sheaves on the total spaces of elliptic fibrations 
we shall need
the notions of ``relative torsion-freeness'' and ``relative 
(semi)stability''. The
first means that the sheaf is flat over the base and torsion-free 
when restricted
to any fibre, while the second means that it is pure-dimensional and the
restriction to any fibre is
$\mu$-(semi)stable in the sense of Simpson \cite{Simp}.

To prove that relatively semistable sheaves of positive degree are
WIT$_0$ we need a preliminary result (cf. also \cite{Br}).

\begin{lemma} A coherent sheaf $\F$ on $X$ is WIT$_0$ if and only if 
$$
\operatorname{Hom}_X(\F, \Pc_\xi)=0
$$
for every $\xi\in \widehat X$.
\label{l:wit}
\end{lemma} 
\begin{proof} By Corollary (\ref{cor1}), $\Pc_\xi$ is WIT$_1$ and $\bS^1(\Pc_\xi)=\kappa(\xi^\ast)$, where $\xi^\ast$ is the point of $\widehat X_t$ corresponding to $\Pc_\xi^\ast$ ($t=\hat p(\xi)$). Then, Parseval theorem  \cite[Prop. 3.4]{BBHM} implies that
$$
\operatorname{Hom}_X(\F, \Pc_\xi)\simeq \operatorname{Hom}_{D(\widehat X)}(\bS(\F), \kappa(\xi^\ast)[1])\,.
$$
If $\F$ is not WIT$_0$, there is a point $\xi^\ast\in\widehat X$ such that a nonzero morphism $\bS^1(\F)\to \kappa(\xi^\ast)$ exists. This gives rise to a non-zero morphism $\bS(\F) \to \kappa(\xi^\ast)[1]$ in the derived category, so that $\operatorname{Hom}_X(\F, \Pc_\xi)\neq 0$. The converse is straightforward.
\end{proof}

\begin{theorem} Let $\F$ be a relatively (semi)stable sheaf on $X$, 
flat over $B$,
with
$d(\F)>0$. Then $\F$ is WIT$_0$ and its Fourier-Mukai transform $\widehat\F$ is
relatively (semi)stable.\label{ss}
\end{theorem}
\begin{proof} We deal with the semistable case, the stable case 
differing only by
some inequalities which become strict. By base change it is enough to restrict to a fibre $X_t$. The fact that $\F$ is WIT$_0$ follows from Lemma \ref{l:wit} since $\F$ is semistable of positive degree.

The sheaf $\widehat \F=\bS^0_t(\F)$ has rank $d(\F)>0$ by Proposition 
\ref{chern}
and is WIT$_1$. All of its subsheaves are WIT$_1$ as well so they cannot have
zero-dimensional support; it follows that
$\widehat \F$ is torsion-free. Let us consider an exact sequence $$ 0\to\E
\to\widehat\F\to\G \to 0
$$ with $\G$ semistable. By applying the inverse Fourier-Mukai transform one
obtains that
$\E$ is WIT$_1$ and that there is an exact sequence $$ 0\to\widehat\bS_t^0\G
\to\widehat\E
\to\F \to\widehat\bS_t^1\G \to 0 $$ which splits into
$$ 0\to\widehat\bS_t^0\G \to\widehat\E \to\Nc\to 0\,,\qquad 0\to\Nc \to \F
\to\widehat\bS_t^1\G \to 0\,.
$$ Then $\widehat\bS_t^1\G$ and $\Nc$ are WIT$_0$, and we obtain 
exact sequences
$$ 0\to\E \to\widehat\Nc \to\bS_t^1\widehat\bS_t^0\G \to 0\,,\qquad 
0\to\widehat\Nc
\to\widehat\F \to\bS_t^0\widehat\bS_t^1\G \to 0\,, $$ so that 
$\widehat\bS_t^0\G$
is WIT$_1$, and one has $$ 0\to\bS_t^1\widehat\bS_t^0\G \to\G
\to\bS_t^0\widehat\bS_t^1\G
\to0\,. $$ If $\bS_t^1\widehat\bS_t^0\G\ne 0$ the semistability of $\G$ implies
that
$d(\bS_t^1\widehat\bS_t^0\G)<0$ which contradicts Corollary \ref{wit} since
$\bS_t^1\widehat\bS_t^0\G$ is WIT$_0$. Then $\bS_t^1\widehat\bS_t^0\G=0$ and we
have
$\E\simeq\widehat\Nc$, so that $\widehat\E\simeq\Nc$ and
$\mu(\widehat\E)\le\mu(\F)$ since
$\F$ is semistable. By Corollary \ref{wit}, $d(\E)\le 0$.

We can exclude the case
$d(\E)=0$ (otherwise
$\widehat\E$ is a subsheaf of $\F$ of rank zero). Finally, if $d(\E)<0$ we can
apply Corollary \ref{mu} to obtain $\mu(\E)\le\mu(\widehat\F)$ which proves the
semistability of
$\widehat\F$. \end{proof}

\begin{corollary} Let $\F$ be a torsion-free semistable sheaf on $X_t$ of degree $d>0$. Then
$H^1(X_t,\F\otimes\Pc_\xi)=0$ for every $\xi\in{\widehat X}_t$.
\end{corollary}
As a side result we show that the relative Fourier-Mukai transform provides a
characterization of some moduli spaces of relatively stable bundles. 
Let $J_n\to
B$ be the relative Jacobian of invertible sheaves of relative degree 
$n$ and $\bar
J_n\to B$ the natural compactification of $J_n$ obtained by adding to 
$J_n$ the $
B$-flat coherent sheaves on $p\colon X\to B$ whose restrictions to 
the fibres of
are torsion-free, of rank one and degree
$n$ \cite{AK}. Theorem \ref{ss} gives:
\begin{corollary} Let $\Nc$ be an invertible sheaf on $X$ of relative 
degree $m$,
and let
$\mathcal M(n,$ $nm-1)$ be the moduli space of  rank $n$ relatively 
$\mu$-stable
sheaves on
$X\to B$ of degree $nm-1$. The Fourier-Mukai  transform induces an 
isomorphism of
$B$-schemes
\begin{eqnarray*}
\bar J_n& \stackrel{\bS^0\otimes\Nc}{\relbar\joinrel\longrightarrow}& 
\mathcal M
(n,nm-1)\\
\L&\mapsto& \bS^0(\L)\otimes\Nc\,.
\end{eqnarray*}
\label{cor3.5}\end{corollary}

\section{Flag schemes and Hecke correspondences\label{hecke}} As a first
application of the above results we construct the relative Hecke 
correspondences
over an elliptic fibration, showing also how, via Fourier-Mukai transform, they
relate to some varieties parametrizing flags of relatively stable sheaves.

\subsection{Flag schemes} We consider exact sequences $$
0\to\L\otimes\A_i^{(t)}\to\E_n
\to\E_{n-i}\to 0 $$ of sheaves on a fibre $X_t$, where $\L$ is torsion-free, of
rank 1 and degree zero, $\A_i^{(t)}$ is the rank $i$ Atiyah sheaf on $X_t$ and
$\E_n$,
$\E_{n-i}$ are torsion-free and stable of degree $d\ge 0$ and rank 
$n$ and $n-i$,
respectively.

We construct a moduli scheme $F^i(n,d)$ parametrizing such exact sequences. The
relevant functor of points is defined as follows: for every 
$B$-scheme $S\to B$,
let ${\mathbf F}^i(n,d)(S)$ be the set of all exact sequences
\begin{equation} 0\to\L\otimes\A_i\to\E_n \to\E_{n-i}\to 0\label{e:flag}
\end{equation} of sheaves on $X_S\to S$, where $\L$ is relatively torsion-free
rank 1 and zero degree, $\A_i$ is the rank $i$ relative Atiyah sheaf 
(Definition
\ref{relatiyah}), and
$\E_n$,
$\E_{n-i}$ are relatively torsion-free and stable sheaves of degree 
$d\ge 0$ and
rank $n$,
$n-i$ respectively. Two such sequences
\begin{eqnarray*} & 0\to\L\otimes\A_i\to \E_n \to\E_{n-i}\to 0 & \\ &
0\to\L'\otimes\A_i\to
\E'_n \to\E'_{n-i}\to 0 & \end{eqnarray*} are considered to be 
equivalent if there
exist line bundles $\L_1$, $\L_2$ on $S$ such that the first sequence is
isomorphic to
$$ 0\to\L'\otimes\A_i\otimes p_S^\ast\L_1 \to\E'_n\otimes p_S^\ast\L_2
\to\E'_{n-i}\to 0
$$ in the obvious sense.

We assume that $n$, $d$ are coprime, so that $\mathcal M(n,d)$ is a fine moduli
scheme (consisting only of stable sheaves) and there is a universal sheaf
$U_{n,d}$ on
$X\times_B\mathcal M(n,d)\to\mathcal M(n,d)$.

\begin{proposition} If $n$ and $d$ are coprime, there is a $B$-scheme 
$F^i(n,d)$
$\to B$ which represents the functor ${\mathbf F}^i(n,d)$.
\label{flags}
\end{proposition}
\begin{proof} Let us write $Y={\widehat X}\times_B\mathcal M(n,d)$, so that
$Y_X=X\times_B{\widehat X}\times_B\mathcal M(n,d)$, and denote by 
$\pi_{jk}$ the
projection of $Y_X$ onto the product of the $j$th and $k$th factors. 
By (1.1) of
\cite{AK} there exists a coherent sheaf ${\mathcal H}$ on $Y$ such 
that for every
morphism
$S\to Y$ and every quasi-coherent sheaf $\F$ on $S$ one has \begin{equation}
\mbox{\rm Hom}_S({\mathcal H}_S,\F)\simeq\mbox{\rm Hom}_{X_S}
((\pi_{12}^\ast\Pc\otimes\pi_1^\ast\A_i)_S, \pi_{13}^\ast
U_{n,d}\otimes_S\F)\,,\label{e:Hsheaf} \end{equation} functorially in $\F$.

If $q\colon P=\mbox{Proj}\,S({\mathcal H})\to Y$ is the associated Grothendieck
projective bundle, there is a universal epimorphism $q^\ast {\mathcal H}
\to\O_P(1)$. According to Eq. (\ref{e:Hsheaf}) this yields a 
universal morphism of
sheaves on $X_P=X\times_B P\to P$
$$ (\pi_{12}^\ast\Pc\otimes\pi_1^\ast\A_i)_P\stackrel{\phi}{\longrightarrow}
\pi_{13}^\ast U_{n,d}\otimes_P\O_P(1)\,.
$$ By (2.3) of \cite{AK} there is a subscheme $V$ of $P$ 
parametrizing the points
$T\to P$ such that
$\mbox{coker}\,\phi_T$ is flat over $T$ and relatively torsion-free 
of rank $n-r$.
The subscheme of $V$ determined by the points where $\mbox{coker}\,\phi$ is
semistable is the desired scheme $F^r(n,d)$.
\end{proof}

We call $F^i(n,d)$ the {\it flag scheme} of order $i$ of $p\colon 
X\to B$. We have
scheme morphisms \begin{equation} {\widehat X}\times_B\M(n,d)
\stackrel{\pi_n}{\longleftarrow} F^i(n,d) \stackrel{\varphi_n}{\longrightarrow}
\M(n-i,d)
\label{e:flagdiag}\end{equation} defined by $\pi_n(0\to\L\otimes\A_i\to\E_n
\to\E_{n-i}\to 0)=(\L,\E_n)$,
$\varphi_n(0\to\L\otimes\A_i\to\E_n \to\E_{n-i}\to 0)=\E_{n-i}$. 

\subsection{Relative Hecke correspondences} Let $\G_r$ be a 
torsion-free sheaf
on $X$ of rank $r$ and relative degree
$m<0$.
\begin{definition} A higher modification of length $i\le r$ at a 
point $x\in X$ is
an exact sequence
\begin{equation} 0\to\G_r\to\G'_r\to\K(x)\to 0\,,\label{e:higmod} 
\end{equation}
where
$\G'_r$ is a torsion-free sheaf of rank $r$ and $\K(x)$ is a coherent torsion
sheaf of length $i$ concentrated at $x$. \end{definition} The sheaf $\G'_r$ has
degree $m-i<0$. We shall construct a $B$-scheme $H^i(r,m)$ parametrizing the
length $i$ higher modifications.

We define now the relevant functor of points. For every morphism
$S\to B$ we let
${\mathbf H}^i(r,m)(S)$ be the family of all exact sequences $$
0\to\G_r\to\G'_r\to\K_i(x)\to 0\,,
$$ of $S$-flat sheaves on $X_S$; here $\G_r$, $\G'_r$ are relatively stable
sheaves of rank
$r$ and negative degree $m$, $m-i$ respectively, while $\K_i(x)$ has 
length $i$ on
every fibre and is concentrated on the image $x(S)$ of a section $x\colon
S\hookrightarrow X_S$ of
$p_S$. Two such sequences
$$\begin{array}{c} 0\to\G_r\to\G'_r\to\K_i(x)\to 0
\\[5pt] 0\to\bar\G_r\to\bar\G'_r\to\bar\K_i(x)\to 0\,,
\end{array}$$ are considered equivalent if there exist line bundles 
$\L_1$, $\L_2$
over $S$ such that
$$ 0\to\G_r\to\G'_r\to\K_i(x)\to 0
$$ is isomorphic to
$$ 0\to\bar\G_r\otimes p_S^\ast\L_2 \to\bar\G'_r \to\bar\K_i(x)\otimes
p_S^\ast\L_1 \to 0\,.
$$

Let us consider integers $d>0$, $n\ge i\ge 0$. \begin{theorem} The 
Fourier-Mukai
transform induces an isomorphism of functors
$$
\bS\colon {\mathbf F}^i(n,d)\simeq {\mathbf H}^i(d,-n)\,. $$
\label{main}\end{theorem}
\begin{proof} The Fourier-Mukai transform interchanges the exact sequences
(\ref{e:flag}) and (\ref{e:higmod}). \end{proof}
\begin{corollary} If $n$, $d$ are coprime, the functor ${\mathbf H}^i(d,-n)$ is
representable by a
$B$-scheme $H^i(d,-n)$, called the length $i$ Hecke correspondence. The
Fourier-Mukai transforms yields a scheme isomorphism
$$
\bS\colon F^i(n,d)\simeq H^i(d,-n)
$$ between the flag scheme and the length $i$ Hecke correspondence.
\label{maincorol}
\end{corollary}

The reason for giving $H^i(d,-n)$ this name is that for every point 
$x\in X$, the
fibre
$H_x^i(d,-n)$ is the graph of the classical length $i$ Hecke correspondence
between the moduli spaces $\M(X_t,d,-n)$ and $\M(X_t,d,i-n)$ on the curve
$X_t$ ($t=p(x)$) (cf. \cite{Dr,La}).

There are scheme morphisms
\begin{equation} X\times_B\M(d,-n) 
\stackrel{\alpha_i}{\longleftarrow} {H^i}(d,-n)
\stackrel{\beta_i}{\longrightarrow}
\M(d,i-n)
\label{e:hecdiag}\end{equation} defined by
$$\begin{array}{c}
\alpha_i(0\to\G_d \to \G'_d
\to\K_i(x)\to 0)=(x,\G_d)\,, \\[7pt]
\beta_i(0\to\G_d \to \G'_d
\to\K_i(x)\to 0) =\G'_d\,.
\end{array}$$

The Fourier-Mukai transform maps isomorphically diagram (\ref{e:flagdiag}) to
diagram (\ref{e:hecdiag}): $$
\xymatrix{ & F^i(n,d) \ar[dr]^{\pi_n} \ar[dl]_{\varphi_n} \\ 
\M(n-i,d) && {\widehat
X}\times_B\M(n,d) }
$$
$$ \stackrel{\bS}{\longrightarrow}
\xymatrix{ & H^i(d,-n) \ar[dr]^{\alpha_i} \ar[dl]_{\beta_i} \\ \M(d,i-n) &&
X\times_B\M(d,-n) }
$$ The isomorphism between ${\widehat X}\times_B\M(n,d)$ and 
$X\times_B\M(d,-n)$
is provided by the morphism
$(\varpi^{-1}\circ\iota,\bS^0)$, where $\iota\colon{\widehat 
X}\to{\widehat X}$ is
the morphism that maps $\xi=[\L]$ into $\xi^\ast=[\L^\ast]$, while 
the isomorphism
between
$\M(n-i,d)$ and $\M(d,i-n)$ is induced by $\bS^0$.

\begin{remark} This result amounts to a generalization to the elliptic relative
case of the geometric Langlands correspondence constructed by Laumon 
\cite{La} in
the case of curves over a field. In Drinfeld's terminology, Theorem \ref{main}
means that the Fourier-Mukai transform induces a duality between the 
operators $S$
and $T_p$ as defined in \cite{Dr}.
\end{remark}

\def\vol#1,{{\bf #1}}
\def\year#1,{(#1),}

\end{document}